\documentclass[11pt, intlimits]{amsart}
\usepackage{amsmath,amssymb}
\usepackage{graphicx}

\begin{document}

\newtheorem{theorem}{Theorem}[section]
\newtheorem{cor}[theorem]{Corollary}
\newtheorem{lem}[theorem]{Lemma}
\newtheorem{lemma}[theorem]{Lemma}
\newtheorem{prop}[theorem]{Proposition}
\newtheorem{proposition}[theorem]{Proposition}
\newtheorem{defn}[theorem]{Definition}
\newtheorem*{remark}{Remark}
\newtheorem{conj}[theorem]{Conjecture}
\newtheorem{ex}[theorem]{Example}

\newcommand{\Z}{{\mathbb Z}} 
\newcommand{\Q}{{\mathbb Q}}
\newcommand{\R}{{\mathbb R}}
\newcommand{\C}{{\mathbb C}}
\newcommand{\N}{{\mathbb N}}
\newcommand{\FF}{{\mathbb F}}
\newcommand{\fq}{\mathbb{F}_q}
\newcommand{\rmk}[1]{\footnote{{\bf Comment:} #1}}

\newcommand{\Vari}{{\operatorname{Var}}}
\newcommand{\Poly}{{\operatorname{Poly}}}
\renewcommand{\mod}{\;\operatorname{mod}}
\newcommand{\ord}{\operatorname{ord}}
\newcommand{\TT}{\mathbb{T}}
\renewcommand{\i}{{\mathrm{i}}}
\renewcommand{\d}{{\mathrm{d}}}
\renewcommand{\^}{\widehat}
\newcommand{\HH}{\mathbb H}
\newcommand{\Vol}{\operatorname{vol}}
\newcommand{\area}{\operatorname{area}}
\newcommand{\tr}{\operatorname{tr}}
\newcommand{\norm}{\mathcal N} 
\newcommand{\intinf}{\int_{-\infty}^\infty}
\newcommand{\ave}[1]{\left\langle#1\right\rangle} 
\newcommand{\Var}{\operatorname{Var}}
\newcommand{\Prob}{\operatorname{Prob}}
\newcommand{\sym}{\operatorname{Sym}}
\newcommand{\disc}{\operatorname{disc}}
\newcommand{\CA}{{\mathcal C}_A}
\newcommand{\cond}{\operatorname{cond}} 
\newcommand{\lcm}{\operatorname{lcm}}
\newcommand{\Kl}{\operatorname{Kl}} 
\newcommand{\leg}[2]{\left( \frac{#1}{#2} \right)}  
\newcommand{\Li}{\operatorname{Li}}
\newcommand{\Sgn}{{\operatorname{Sgn}}}
\newcommand{\uuu}{{\operatorname{o}}}
\newcommand{\sumstar}{\sideset \and^{*} \to \sum}

\newcommand{\LL}{\mathcal L} 
\newcommand{\sumf}{\sum^\flat}
\newcommand{\Hgev}{\mathcal H_{2g+2,q}}
\newcommand{\USp}{\operatorname{USp}}
\newcommand{\conv}{*}
\newcommand{\dist} {\operatorname{dist}}
\newcommand{\CF}{c_0} 
\newcommand{\kerp}{\mathcal K}

\newcommand{\Cov}{\operatorname{cov}}
\newcommand{\Sym}{\operatorname{Sym}}

\newcommand{\ES}{\mathcal S} 
\newcommand{\EN}{\mathcal N} 
\newcommand{\EM}{\mathcal M} 
\newcommand{\Sc}{\operatorname{Sc}} 
\newcommand{\Ht}{\operatorname{Ht}}

\newcommand{\aap}{\operatorname{primes}}
\newcommand{\E}{\operatorname{E}} 
\newcommand{\sign}{\operatorname{sign}} 

\newcommand{\divid}{d} 
\newcommand{\inv}{\theta}

\newcommand{\diag}{\operatorname{diag}}
\newcommand{\monic}{\operatorname{monic}}

\newcommand{\acks}{\textit{Acknowledgements}}

\title{Arithmetic correlations over large finite fields  }
\date{\today}

\author{J.P. Keating and E. Roditty-Gershon}

\address{School of Mathematics, University of Bristol, Bristol BS8 1TW, UK}
\email{j.p.keating@bristol.ac.uk}

\address{School of Mathematics, University of Bristol, Bristol BS8 1TW, UK}
\email{(Correspondence to be sent to) er14265@bristol.ac.uk}


\maketitle
\setcounter{page}{1}

\begin{abstract}
The auto-correlations of arithmetic functions, such as the von Mangoldt function, the M$\ddot{\uuu}$bius function and the divisor function, are the subject of classical problems in analytic number theory. The function field analogues of these problems have recently been resolved in the limit of large finite field size $q$. However, in this limit the correlations disappear: the arithmetic functions become uncorrelated. We compute averages of terms of lower order in $q$ which detect correlations.  Our results show that there is considerable cancellation in the averaging and have implications for the rate at which correlations disappear when $q \rightarrow\infty$; in particular one cannot expect remainder terms that are of the order of the square-root of the main term in this context.
\end{abstract}

\section{Introduction}

The generalized twin-prime conjecture \cite{hardy}, predicts that for an $\textit{r}$-tuple of distinct integers $a_{1},\ldots,a_{r}$, which do not cover all residues modulo some prime $p$, there are infinitely many positive integers $n$ such that $n+a_{i}$ are primes, for all $1\leq i \leq r.$ In other words, for $a=(a_{1},\ldots,a_{r})$ let
\begin{equation}
\pi(x;a)=\# \{1\leq n \leq x |n+a_{1},\ldots,n+a_{r} \aap \}
\end{equation}
then the conjecture says that $\pi(x;a)\rightarrow \infty$ as $x\rightarrow\infty,$ unless the local obstruction noted above holds.
If $a_{1},\ldots,a_{r}$ cover all residues modulo some prime $p$, then for any $n$ there is an $i$ such that $p|n+a_{i}$ and $\pi(x;a)$ is bounded as $x\rightarrow\infty.$ Note that the number of primes up to $x,$  $\pi(x),$ is equal to $\pi(x;0),$ and the number of twin primes is $\pi(x;0,2).$

The Hardy-Littlewood conjecture \cite{hardy}, gives the rate in which $\pi(x;a)$ tends to infinity: let
\begin{equation}
\Upsilon(p;a)=\# \{a_{1} \mod p,\ldots, a_{r} \mod p\}
\end{equation}
and
\begin{equation}
C_{r}(a)=\prod_{p}\frac{1-\Upsilon(p;a)/p}{(1-1/p)^{r}}
\end{equation}
then, unless  $a_{1},\ldots,a_{r}$ cover all residues modulo some prime $p,$ the product converges, i.e. $C_{r}(a)>0$, and the conjecture predicts that 
\begin{equation}
\pi(x;a)\sim C_{r}(a)\frac{x}{\log^{r}x},~~~~~~~~~x\rightarrow\infty
\end{equation}  
The Hardy-Littlewood conjecture can be rephrased using the von Mangoldt function
\begin{equation}\label{von}
\Lambda(n) = \begin{cases} \log p & \mbox{if }n=p^k \mbox{ for some prime } p \mbox{ and integer } k \ge 1, \\ 0 & \mbox{otherwise} \end{cases}
\end{equation}
as follows:
\begin{equation}
\sum_{n\leq x}\Lambda(n+a_{1})\cdots\Lambda(n+a_{r})\sim  C_{r}(a)x
\end{equation}

There has recently been interest in the analogue of the Hardy-Littlewood conjecture over large finite fields.
Let $ \fq $ be a finite field of $q$ elements, with $q$ odd, and let $\fq[T] $ be the ring of polynomials with coefficients in $ \fq $. Let $ P_{n}:=\{f\in \fq[T]:\deg f=n\} $ be the set of polynomials of degree $  n$, and $M_{n}:=\{f\in \fq[T]:\deg f=n, f \monic\}$ be the set of monic polynomials of degree $  n.$

It follows from the work of Bary-Soroker \cite{lior} 
and Bender and Pollack \cite{pollack}, that for fixed $n$ and in the limit of large (odd) $q$, the analogue of the Hardy-Littlewood conjecture holds in the form
\begin{equation}\label{hardy-littlewood}
\frac{1}{q^{n}}\sum_{f\in M_{n}} \Lambda(f)\Lambda(f+K)=1+E(K,n,q)
\end{equation}
with
\begin{equation}\label{error}
E(K,n,q)=O(q^{-\frac{1}{2}})
\end{equation}
and where, by analogy with \eqref{von}
\begin{equation}
\Lambda(f) = \begin{cases} \deg p & \mbox{if }f=p^k \mbox{ for some prime polynomial } p \mbox{ and integer } k \ge 1, \\ 0 & \mbox{otherwise} \end{cases}
\end{equation}
and $0\neq K\in \fq[T]$, $n>\deg K$. The bound \eqref{error} results from using an algebraic irreducibility criterion and then invoking the Lang-Weil bound for the number of points on varieties in finite fields. 
It is significant that at leading order as $q\rightarrow\infty$, the von Mangoldt functions are uncorrelated in that the autocorrelation function is independent of $K$. Information about any correlations is therefore contained in the error term $ E(K,n,q);$ however the bound \eqref{error} is not sensitive enough to detect this. 

It is natural to speculate that the true order of $E$ corresponds to square-root cancellation in the sum in \eqref{hardy-littlewood}, so that $E(K,n,q)=O(q^{-n/2})$.  Our results imply that this cannot, in fact, be the case.

In this note, we study the sum of the error term $ E(K,n,q) $ over all monic polynomials $K$ of a given degree.
The method we use here is based on a calculation of Keating and Rudnick \cite{keating-rudnick 1},
who computed the variance of the sum of the von Mangoldt function over short intervals and in arithmetic progressions. Our approach applies as well to other arithmetic functions such as the 
M$\ddot{\uuu}$bius function $\mu(f)$ and the divisor function $d(f)$. We begin by reviewing the result of \cite{keating-rudnick 1}, which we shall need here. We then state our results for $\Lambda(f)$. The extensions to $\mu(f)$ and $d(f)$ are outlined in section 4.
\subsection{Arithmetic progressions}
For a polynomial $Q\in \fq[T]$ of a positive degree, and $A\in \fq[T]$ coprime to $Q$ and any $n>0$, set 
\begin{equation}
\Psi(n;Q,A):=\sum_{\substack{N\in M_{n}\\ N=A \mod Q}}\Lambda(N)
\end{equation}
The prime polynomial theorem in arithmetic progressions states that as $n\rightarrow \infty$
\begin{equation}\label{pnt in arithmetic}
\Psi(n;Q,A)\sim  \frac{q^{n}}{\Phi(Q)}
\end{equation}
where $\Phi(Q)$ is the Euler totient function for this context, namely the number of reduced residue classes modulo $Q.$ Now set
\begin{equation}
G(n;Q):=\sum_{\substack{A \mod Q \\ \gcd(A,Q)=1}}|\Psi(n;Q,A)- \frac{q^{n}}{\Phi(Q)}|^{2}
\end{equation}
It was shown in \cite{keating-rudnick 1} that the following holds
\begin{theorem}\label{k-r progressions}
In the limit $q\rightarrow \infty,$ 
\begin{equation}\label{main k-r prog}
\frac{G(n;Q)}{q^{n}}=(\deg Q-1)+O(\frac{1}{\sqrt{q}})
\end{equation}
where $\deg Q \leq n.$ 
\end{theorem}
The result is based on an equidistribution statement for the unitarized Frobenii of primitive odd characters \cite{katz 1}.
\subsection{Short intervals}
For $A\in P_{n}$ of degree $n$, and $h<n$, a short interval of size $h$ around $A$ is defined by
\begin{equation}
I(A;h):=\{f: ||f-A||\leq q^{h}\}=A+P_{\leq h} 
\end{equation}
Where $||f||:=q^{\deg f}$ and $P_{\leq h}=\{0\}\cup \bigcup_{0\leq m \leq h}P_{m} $ is the space of polynomials of degree at most $h$. We have 
\begin{equation}
\#I(A;h)=q^{h+1}
\end{equation}
Note that if $f, g \in I(A;h)$, then there exists a polynomial of degree less than or equal to $h$ such that $f$ and $g$ are congruent modulo this polynomial.
For $1\leq h <n$ and $A\in P_{n}$, set
\begin{equation}
 \upsilon(A;h)=\sum_{\substack{f\in I(A;h)}}\Lambda(f)
 \end{equation} 
The mean value of $\upsilon(A;h)$, when we average over $M_{n},$ is 
\begin{equation}
\begin{split}
\langle\upsilon(\bullet;h)\rangle&:=\frac{1}{q^{n}}\sum_{A\in M_{n}}\upsilon(A;h)
\\&=\frac{1}{q^{n}}\sum_{A\in M_{n}}\sum_{\substack{f\in I(A;h)}}\Lambda(f)
\\&=\frac{1}{q^{n}}q^{h+1}\sum_{f\in M_{n}}\Lambda(f)
\\&=q^{h+1}
\end{split}
\end{equation}
The last equality is due to the Prime Polynomial Theorem, which in this context is the identity
\begin{equation}
\sum_{f\in M_{n}}\Lambda(f)=q^{n}
\end{equation}
Now, consider the limit as $q\rightarrow \infty$ of the variance of $\upsilon(A;h)$
\begin{equation}
V(\upsilon(\bullet;h))=\frac{1}{q^{n}}\sum_{A\in M_{n}}|\upsilon(A;h)-\langle\upsilon(\bullet;h)\rangle|^{2}
\end{equation}
The following theorem was established in \cite{keating-rudnick 1} 
\begin{theorem}\label{k-r short intervals}
In the limit $q\rightarrow \infty$, and for $h<n-3$, 
\begin{equation}\label{k-r}
\frac{1}{q^{h+1}}V(\upsilon(\bullet;h))=n-h-2+O(\frac{1}{\sqrt{q}})
\end{equation}
\end{theorem}
The proof is based on an equidistribution statement for the unitarized Frobenii of primitive even characters \cite{katz 2}. 

\subsection{Statement of results}
The first quantity we study in this note is the sum of the error term $E(K,n,q)$, defined in \eqref{error}, over all monic polynomials of a given degree: 
\begin{equation}\label{reg sum}
S_{E}(k,n,q):=\sum_{K\in M_{k}}E(K,n,q)
\end{equation}
The second quantity we study is a "twisted" sum of $E(K,n,q)$: 
\begin{equation}\label{twisted sum}
\tilde{S}_{E}(n,q;Q)=\sum_{j=0}^{n-\deg Q-1}\sum_{K\in M_{j}}E(KQ,n,q)
\end{equation}
We have the following theorems for the above sums:
\begin{theorem}\label{main1}
For $k<n-3$ and in the limit $q\rightarrow\infty$, 
\begin{equation}\label{main1eq}
S_{E}(k,n,q)=\frac{1}{1-q}+O(\frac{1}{q^{3/2}})
\end{equation}
\end{theorem}

\begin{theorem}\label{main2}
Let $Q$ be a polynomial such that $\deg Q<n$. Then in the limit $q\rightarrow\infty$, 
\begin{equation}
\tilde{S}_{E}(n,q;Q)=\frac{(n-\deg Q)}{1-q}+O(\frac{1}{q^{3/2}})
\end{equation}
\end{theorem}

\begin{cor}\label{cor main2}
For a polynomial $Q$ of degree $n-1$ the error term is
\begin{equation}
E(Q,n,q)=-\frac{1}{q-1}+O(\frac{1}{q^{3/2}})
\end{equation}
\end{cor}
We note that the number of terms in the sum in \eqref{reg sum} is $q^{k}$, so the average of $E(K,n,q)$ over $K\in M_{k}$ is approximately $-\frac{1}{q^{k}(q-1)}.$ Similarly, the number of terms in the sum in \eqref{twisted sum} is $\frac{q^{n-k}-1}{q-1}$ where $k=\deg Q$, and so the average of $E(KQ,n,q)$ over $K\in M_{<n-k}$ is approximately $\frac{n-k}{q^{n-k}-1}.$

Theorems \ref{main1} and \ref{main2} have implications for $E(K, n, q)$ that are worth noting here.  First, they both involve sums of a large number of terms ($q^k$ and $O(q^{n-k-1})$ respectively).  Therefore, the fact that both vanish as $q \rightarrow \infty$ necessitates a considerable degree of cancellation.  Second, setting $k< n/2-1$, one sees that \eqref{main1eq} is not consistent with $E(K, n, q)$ being $O(q^{-n/2})$, as might have been anticipated based on the (optimistic) assumption of square-root cancellation in the fluctuations in the sum in \eqref{hardy-littlewood}. Indeed, setting $k=0$, one sees that $E(1, n, q)=\frac{1}{1-q}+O(\frac{1}{q^{3/2}})$ for all $n>4$. The same implication follows from corollary \ref{cor main2}.  More generally, the maximum of $|E(K, n, q)|$ with respect to $K\in M_k$ is bounded from below by the modulus of the average of $E$, which is  $\frac{1}{q^{k}(q-1)}$, and this is larger than $q^{-n/2}$ for $k<n/2-1$.

\section{Sum of error terms}
In this section we use the result obtained in Theorem \ref{k-r short intervals}, to evaluate the sum of the error terms $ E(K,n,q).$

First, we use the definition of the variance of $\upsilon(A;h)$ to write
\begin{equation}
V(\upsilon(\bullet;h))=\frac{1}{q^{n}}[\sum_{A\in M_{n}}\upsilon(A;h)^{2}-2q^{h+1}\sum_{A\in M_{n}}\upsilon(A;h)+q^{n}q^{2(h+1)}]
\end{equation}
Note that
\begin{equation}
\begin{split}
\sum_{A\in M_{n}}\upsilon(A;h)=q^{n}\langle\upsilon(\bullet;h)\rangle=q^{h+1+n}
\end{split}
\end{equation}
For the sum involving $\upsilon(A;h)^{2}$ we have
\begin{equation}
\begin{split}
\sum_{A\in M_{n}}\upsilon(A;h)^{2}&=\sum_{A\in M_{n}}\sum_{\substack{f,g\in I(A;h)}}\Lambda(f)\Lambda(g)
\\&=\sum_{A\in M_{n}}\sum_{\substack{f\in I(A;h)}}\Lambda(f)^{2}+\sum_{A\in M_{n}}\sum_{\substack{f,g\in I(A;h)\\f\neq g}}\Lambda(f)\Lambda(g)
\end{split}
\end{equation}

For the first term, we have (see Lemma 3.1 in \cite{keating-rudnick 1})
\begin{equation}
\begin{split}
\sum_{A\in M_{n}}\sum_{\substack{f\in I(A;h)}}\Lambda(f)^{2}&=
q^{h+1}\sum_{f\in M_{n}}\Lambda(f)^{2}
\\&=q^{h+n+1}n+O(n^{2}q^{n/2})
\end{split}
\end{equation}

For the second term, recall that if $f, g \in I(A;h)$, then there exists a polynomial of degree smaller or equal to $h$ such that $f$ and $g$ are congruent modulo this polynomial. The sum over $f\neq g$ thus can be written as
\begin{equation}
\sum_{A\in M_{n}}\sum_{\substack{f\in I(A;h)}}\sum_{j=0}^{ h}\sum_{\substack{\deg J=j\\J\neq 0}}\Lambda(f)\Lambda(f+J)=
q^{h+1}\sum_{\substack{f\in M_{n}}}\sum_{j=0}^{ h}\sum_{\substack{\deg J=j\\J\neq 0}}\Lambda(f)\Lambda(f+J)
\end{equation}
We will restrict the $J$-sum to monics, multiplying it by $q-1$, and so the right-hand side becomes
\begin{equation}\label{monic}
\begin{split}
&q^{h+1}(q-1)\sum_{\substack{f\in M_{n}}}\sum_{j=0}^{ h}\sum_{\substack{ J\in M_{j}}}\Lambda(f)\Lambda(f+J)
\end{split}
\end{equation}
Combining the above, we get
\begin{equation}
\begin{split}
\frac{V(\upsilon(\bullet;h))}{q^{h+1}}&=n-q^{h+1}+(q-1)\frac{1}{q^{n}}\sum_{f\in M_{n}}\sum_{j=0}^{h}\sum_{ J\in M_{j}}\Lambda(f)\Lambda(f+J)+O(\frac{1}{q^{n/2+h+1}})
\end{split}
\end{equation}

By the definition $\eqref{hardy-littlewood}$ of $E(J,n,q)$, we have
\begin{equation}
\begin{split}
\frac{V(\upsilon(\bullet;h))}{q^{h+1}}&=n-q^{h+1}+(q-1)\sum_{j=0}^{h}\sum_{ J\in M_{j}}(1+E(J,n,q))+O(\frac{1}{q^{n/2+h+1}})
\\&=n-1+(q-1)\sum_{j=0}^{ h}\sum_{\substack{ J\in M_{j}}}E(J,n,q)+O(\frac{1}{q^{n/2+h+1}})
\end{split}
\end{equation}

Now, we combine the expansion of $V(\upsilon(\bullet;h))$ with Theorem \ref{k-r short intervals}, to write
\begin{equation}
(q-1)\sum_{j=0}^{ h}\sum_{ J\in M_{j}}E(J,n,q)=-(h+1)+O(\frac{1}{\sqrt{q}})
\end{equation}
By subtracting the sum up to $h-1$ from the sum up to $h$, we get Theorem $\ref{main1}$
\begin{equation}
(q-1)\sum_{ J\in M_{h}}E(J,n,q)=-1+O(\frac{1}{\sqrt{q}})
\end{equation}
Thus the average of the error terms is 
\begin{equation}
\frac{1}{q^{h}}\sum_{ J\in M_{h}}E(J,n,q)=-\frac{1}{(q-1)q^{h}}+O(\frac{1}{q^{h+3/2}})
\end{equation}
\section{Twisted sum}
In this section, we use Theorem \ref{k-r progressions} to evaluate the sum of $E(KQ,n,q)$ for fixed $Q$ with $K\in M_{<n-\deg Q}.$
First, we use the definition of $ G(n;Q) $ to write
\begin{equation}
G(n;Q)=\sum_{\gcd(A,Q)=1}\Psi(n;Q,A)^{2}-2\frac{q^{n}}{\Phi(Q)}\sum_{\gcd(A,Q)=1}\Psi(n;Q,A)+
\frac{q^{2n}}{\Phi(Q)}
\end{equation}
Note that 

\begin{equation}
\begin{split}
\sum_{\gcd(A,Q)=1}\Psi(n;Q,A)&=\sum_{\substack{\deg f =n\\ \gcd(f,Q)=1}}\Lambda(f)\\&=
\sum_{\deg f =n}\Lambda(f)-\sum_{\substack{\deg f =n\\ \deg \gcd(f,Q)>0}}\Lambda(f)\\&=
q^{n}-\sum_{\substack{\deg P |n\\ P|Q \\ prime}}\deg P\\&=
q^{n}+O(\deg Q)
\end{split}
\end{equation}

For the sum over $\Psi(n;Q,A)^{2}$, we have
\begin{equation}
\begin{split}
\sum_{\gcd(A,Q)=1}\Psi(n;Q,A)^{2}&=\sum_{\substack{f,g\in M_{n} \\ f=g\mod Q\\ \gcd(f,Q)=1}}\Lambda(f)\Lambda(g)
\\&=\sum_{\substack{f\in M_{n} \\  \gcd(f,Q)=1}}\Lambda(f)^{2}+\sum_{\substack{f,g\in M_{n} \\ f=g\mod Q\\ \gcd(f,Q)=1 \\ f\neq g}}\Lambda(f)\Lambda(g)
\end{split}
\end{equation}
For the first term, we have (by Lemma 3.1 in \cite{keating-rudnick 1})

\begin{equation}
\begin{split}
\sum_{\substack{f\in M_{n} \\  \gcd(f,Q)=1}}\Lambda(f)^{2}&=\sum_{f\in M_{n} }\Lambda(f)^{2}-\sum_{\substack{\deg f =n\\ \deg \gcd(f,Q)>0}}\Lambda(f)^{2}\\&=
nq^{n}+O(n^{2}q^{n/2})+O(\deg Q ^2)
\end{split}
\end{equation}

The sum over $f\neq g$ can be written as in the appendix of \cite{keating-rudnick 1} (equation A.14)
\begin{equation}
\sum_{\substack{f,g\in M_{n} \\ f=g\mod Q\\ \gcd(f,Q)=1 \\ f\neq g}}\Lambda(f)\Lambda(g)=
(q-1)\sum_{j=0}^{n-\deg Q -1}\sum_{J\in M_{j}}\sum_{f\in M_{n}}\Lambda(f)\Lambda(f+JQ)
\end{equation}
Combining the above, we obtain
\begin{equation}
\begin{split}
\frac{G(n;Q)}{q^{n}}&=\frac{1}{q^{n}}(q-1)\sum_{j=0}^{n-\deg Q -1}\sum_{J\in M_{j}}\sum_{f\in M_{n}}\Lambda(f)\Lambda(f+JQ)
\\&+n+O(n^{2}q^{-n/2})+O(\frac{\deg Q ^2}{q^{n}})-\frac{q^{n}}{\Phi(Q)}+O(\frac{\deg Q}{\Phi(Q)})
\end{split}
\end{equation}
By using the function field version of the Hardy-Littlewood conjecture, i.e. $\eqref{hardy-littlewood}$, and noting that as $q\rightarrow \infty$
\begin{equation}
\frac{q^{\deg Q}}{\Phi(Q)} \rightarrow 1 
\end{equation}
we have
\begin{equation}
\begin{split}
\frac{G(n;Q)}{q^{n}}&=(q-1)\sum_{j=0}^{n-\deg Q -1}\sum_{ J\in M_{j}}(1+E(JQ,n,q))
\\&+n+O(n^{2}q^{-n/2})+O(\frac{\deg Q ^2}{q^{n}})-q^{n-\deg Q}+O(\frac{\deg Q}{\Phi(Q)})\\&=
(q-1)\sum_{j=0}^{n-\deg Q -1}\sum_{ J\in M_{j}}E(JQ,n,q) 
\\&+n-1+O(n^{2}q^{-n/2})+O(\frac{\deg Q ^2}{q^{n}})+O(\frac{\deg Q}{\Phi(Q)})
\end{split}
\end{equation}

Now, we combine the expansion of $G(n;Q)$ with $\eqref{main k-r prog}$, obtaining
\begin{equation}
\begin{split}
(q-1)\sum_{j=0}^{n-\deg Q -1}\sum_{J\in M_{j}}E(JQ,n,q)+n-1=(\deg Q-1)+O(\frac{1}{\sqrt{q}})
\end{split}
\end{equation}
Therefore we have the following expression for the error term $E(JQ,n,q)$:
\begin{equation}\label{twistedfinal}
(q-1)\sum_{j=0}^{n-\deg Q -1}\sum_{J\in M_{j}}E(JQ,n,q)=-(n-\deg Q)+O(\frac{1}{\sqrt{q}})
\end{equation}
which proves Theorem $\eqref{main2}.$
From this we can deduce that for a polynomial $Q$ of degree $n-1$ the error term is
\begin{equation}
(q-1)E(Q,n,q)=-1+O(\frac{1}{\sqrt{q}})
\end{equation}
and that in general the average is $\frac{n-\deg Q}{q^{n-\deg Q}-1}+O(\frac{1}{q^{n-\deg Q+1/2}})$, because the number of terms in the sum in \eqref{twistedfinal} is $\frac{q^{n-\deg Q}-1}{q-1}$.

\section{Further examples -- the M$\ddot{\uuu}$bius function and the divisor function}
In the following section we use the method demonstrated above to evaluate sums of the error terms in two other important problems. The first is the additive divisor problem over $\fq[T]$  (see \cite{Andrade}), and the second is Chowla's conjecture over $\fq[T]$ (see \cite{carmon}).
\subsection{Definitions}
For an arithmetic function $\alpha$ 
we define
\begin{equation}
\upsilon_{\alpha}(A;h)=\sum_{f\in I(A;h)}\alpha(f)
\end{equation}
The mean value of $\upsilon_{\alpha}(A;h)$ is therefore 
\begin{equation}
\begin{split}
\langle\upsilon_{\alpha}(\bullet;h)\rangle &=\frac{1}{q^{n}}\sum_{A\in M_{n}}\upsilon_{\alpha}(A;h) 
\\&=\frac{1}{q^{n}}q^{h+1}\sum_{f\in M_{n}}\alpha(f)
\\&=q^{h+1}\langle\ \alpha \rangle_{n}
\end{split}
\end{equation}
The variance of $\upsilon_{\alpha}(A;h)$ is given by
\begin{equation}
V(\upsilon_{\alpha}(\bullet;h)):=\frac{1}{q^{n}}\sum_{A\in M_{n}}|\upsilon_{\alpha}(A;h)-\langle\upsilon_{\alpha}(\bullet;h)\rangle|^{2}
\end{equation}
In the following we use the function field zeta function
\begin{equation}
\zeta_{q}(s)=\frac{1}{1-q^{1-s}}
\end{equation}
which we also write as $Z(u)=(1-qu)^{-1}$ by setting $u=q^{-s}.$
\subsection{The divisor function}
The additive divisor problem over $\mathbf{Z}$ (sometimes called "shifted divisor" or "shifted convolution" problem), concerns the asymptotics of the sum
\begin{equation}
  D(x;h):=\sum_{n\leq x}d(n)d(n+h)
   \end{equation} 
   where $d$ is the divisor function.
Ingham \cite{Ingham} 
computed the leading term, and Estermann \cite{Estermann} 
gave an asymptotic expansion
\begin{equation}
\sum_{n\leq x}d(n)d(n+h)=xP_{2}(\log x;h)+O(x^{\frac{11}{12}}(\log x)^{3})
\end{equation}
where 
\begin{equation}
P_{2}(u;h)=\frac{1}{\zeta(2)}\sigma_{-1}(h)u^{2}+a_{1}(h)u+a_{2}(h)
\end{equation}
with
\begin{equation}
\sigma_{w}(h)=\sum_{k|h}k^{w}
\end{equation}
and $a_{1}(h),a_{2}(h)$ are complicated coefficients.

The size of the reminder term plays an important role in various problems in analytic number theory; see, for example,  \cite{Iwaniec}, \cite{heath}.

Andrade, Bary-Soroker and Rudnick \cite{Andrade} 
studied the  additive divisor problem over $\fq[T]$, showing that in the limit $q\rightarrow \infty$
\begin{equation}\label{abr}
\frac{1}{q^{n}}\sum_{f\in M_{n}}d(f)d(f+J)=(n+1)^{2}+E_{d}(J,n,q)
\end{equation}
when $0\neq J\in \fq[T]$, and $\deg(J)<n,$ with
\begin{equation}
E_{d}(J,n,q)=O(q^{-\frac{1}{2}})
\end{equation}
This corresponds to $\alpha (f)=d(f)$ in the definition above.
As before, we will use a result that is based on an equidistribution statement for the unitarized Frobenii of primitive even characters in order to study the sum of $E(J,n,q)$ over monic polynomials  of a given degree. To this end, we quote the following theorem from \cite{roditty}:
\begin{theorem}\label{divisor var}
In the limit of large field size, $q\rightarrow\infty,$ the following holds:
If $0\leq h \leq \frac{n}{2}-2$ then
\begin{equation}
V(\upsilon_{d}(\bullet;h))=q^{h+1}\binom{n-2h+1}{3}+O(\frac{q^{h+1}}{\sqrt{q}})
\end{equation}
If $h=\frac{n}{2}-1$ then
\begin{equation}
V(\upsilon_{d}(\bullet;h))=O(\frac{q^{h+1}}{\sqrt{q}})
\end{equation}
If $\frac{n}{2}\leq h < n$ then
\begin{equation}
V(\upsilon_{d}(\bullet;h))=0
\end{equation}
\end{theorem}

We use the definition of the variance of $V(\upsilon_{d}(\bullet;h))$ to obtain
\begin{equation}
\begin{split}
V(\upsilon_{d}(\bullet;h))=\frac{1}{q^{n}}[&q^{h+1}\sum_{f\in M_{n}}d(f)^{2}+q^{h+1}\sum_{f\in M_{n}}\sum_{j=0}^{h}\sum_{\substack{\deg J=j\\J\neq 0}}d(f)d(f+J)
\\&-2q^{2(h+1)}q^{n}(\langle\ d \rangle_{n})^{2}+q^{2(h+1)}(\langle\ d \rangle_{n})^{2}]
\end{split}
\end{equation}
By considering the generating functions of $\sum_{f\in M_{n}}d(f)$ and of $\sum_{f\in M_{n}}d(f)^{2}$, which are $Z(u)^{2}$ and $\frac{Z(u)^{4}}{Z(u^{2})}$ respectively, we have that
\begin{equation}
\sum_{f\in M_{n}}d(f)=q^{n}(n+1)
\end{equation}
and
\begin{equation}
\sum_{f\in M_{n}}d(f)^{2}=q^{n}\binom{n+3}{3}-q^{n-1}\binom{n+1}{3}
\end{equation}
Thus the variance of $V(\upsilon_{d}(\bullet;h))$ is
\begin{equation}
q^{h+1}[\binom{n+3}{3}-q^{-1}\binom{n+1}{3}-q^{h+1}(n+1)^{2}+\frac{1}{q^{n}}\sum_{f\in M_{n}}\sum_{j=0}^{h}\sum_{\substack{\deg J=j\\J\neq 0}}d(f)d(f+J)]
\end{equation}
which by $\eqref{abr}$ is
\begin{equation}
q^{h+1}[\binom{n+3}{3}-q^{-1}\binom{n+1}{3}-(n+1)^{2}+(q-1)\sum_{j=0}^{h}\sum_{ J\in M_{j}}E_{d}(J,n,q)]
\end{equation}
By combining the above with Theorem $\eqref{divisor var}$, and subtracting the sum up to $h-1$ from the sum up to $h$, we have
\begin{theorem}\label{Jsum}
In the limit $q\rightarrow\infty:$ 
\\If $0\leq h \leq \frac{n}{2}-2$ then
\begin{equation}
\sum_{ J\in M_{h}}E_{d}(J,n,q)=\binom{n-2h-1}{3}-\binom{n-2h+1}{3}+O(\frac{1}{\sqrt{q}})
\end{equation}
If $h=\frac{n}{2}-1$ then
\begin{equation}
\sum_{ J\in M_{h}}E_{d}(J,n,q)=-1+O(\frac{1}{\sqrt{q}})
\end{equation}
If $\frac{n}{2}<h<n$ then
\begin{equation}
\sum_{J\in M_{h}}E_{d}(J,n,q)=0
\end{equation}
\end{theorem} 
We note that the number of terms in the sum over $J\in M_{h}$ is $q^{h}$, and so Theorem \ref{Jsum} determines the average of $E_{d}$ when both sides of the equation are divided by $q^h$.  As for the von Mangoldt function, this theorem demonstrates considerable cancellation when $E_{d}(J,n,q)$ is summed over $J$, and establishes a lower bound on its size which rules out the optimistic guess that $E_{d}(J,n,q)=O(q^{-n/2})$.

\subsection{The M$\ddot{\uuu}$bius function}
Chowla's conjecture \cite{Chowla} asserts that given an $\textit{r}$-tuple of distinct integers $\alpha_{1},\ldots,\alpha_{n}$, and $\epsilon_{i}\in\{1,2\},$ not all even, then

\begin{equation}
 \lim_{N\rightarrow\infty}\frac{1}{N}\sum_{n\leq N}\mu(n+\alpha_{1})^{\epsilon_{1}}\cdots\mu(n+\alpha_{r})^{\epsilon_{r}}=0
 \end{equation}                                                                                                                  
This conjecture has recently been shown by Sarnak to imply that $\mu(n)$ does not correlate with any zero entropy sequence \cite{sarnak}. 

Carmon and Rudnick \cite{carmon}, proved the function field version of the Chowla conjecture, in the limit of large (odd) field size. The extension to even characteristics has been established by Carmon in \cite{carmon even}. Here the Chowla conjecture, is shown to hold in the form
\begin{equation}\label{cr}
|\sum_{f\in M_{n}}\mu(f+\alpha_{1})^{\epsilon_{1}}\cdots\mu(f+\alpha_{r})^{\epsilon_{r}}|\leq 2rnq^{n-1/2}+3rn^{2}q^{n-1}
 \end{equation} 
where $r>1,n>1$ (in the case of even characteristics $n>2$) and $\alpha_{1},\ldots,\alpha_{n}\in \fq[T]$ are distinct polynomials with $\deg \alpha_{j}<n.$ As before, $\epsilon_{i}\in\{1,2\},$ not all even.

We will focus on the case of $r=2, \alpha_{1}=0,$ and $\epsilon_{1},\epsilon_{2}=1.$
Denote by 
\begin{equation}
E_{\mu}(J,n,q):=\frac{1}{q^{n}}\sum_{f\in M_{n}}\mu(f)\mu(f+J)
\end{equation}
As before, we will use a result that is based on an equidistribution statement for the unitarized Frobenii of primitive even characters, in order to study the sum of $E_{\mu}(J,n,q)$ over monic polynomials  of a given degree. To this end, we quote the following theorem from \cite{keating-rudnick 2}
\begin{theorem}\label{var mubius}
In the limit of large field size, $q\rightarrow\infty,$ and for $h\leq n-4$,
\begin{equation}
V(\upsilon_{\mu}(\bullet;h))=q^{h+1}+O(q^{h+1/2})
\end{equation}
\end{theorem}

Next, we use the definition of the variance of $V(\upsilon_{\mu}(\bullet;h))$ to obtain
\begin{equation}
\begin{split}
V(\upsilon_{\mu}(\bullet;h))=\frac{1}{q^{n}}[&q^{h+1}\sum_{f\in M_{n}}\mu(f)^{2}+q^{h+1}\sum_{f\in M_{n}}\sum_{j=0}^{h}\sum_{\deg J=j}\mu(f)\mu(f+J)
\\&-2q^{2(h+1)}q^{n}(\langle\ \mu \rangle_{n})^{2}+q^{2(h+1)}(\langle\ \mu \rangle_{n})^{2}]
\end{split}
\end{equation}
The analysis in this case follows exactly the same lines as in the previous calculations and so we omit the details.  
By considering the generating functions of $\sum_{f\in M_{n}}\mu(f)$ and of $\sum_{f\in M_{n}}\mu(f)^{2}$, which are $\frac{1}{Z(u)}$ and $\frac{Z(u)}{Z(2u)}$ respectively, we conclude that for $n>1$
\begin{equation}
\sum_{f\in M_{n}}\mu(f)=0
\end{equation}
\begin{equation}
\sum_{f\in M_{n}}\mu(f)^{2}=\frac{q^{n}}{\zeta_{q}(2)}
\end{equation}
Thus the variance of $V(\upsilon_{d}(\bullet;h))$ is
\begin{equation}
q^{h+1}[\frac{1}{\zeta_{q}(2)}+(q-1)\sum_{j=0}^{h}\sum_{\deg J=j}E_{\mu}(J,n,q)]
\end{equation}
By combining the above with Theorem $\eqref{var mubius}$, and subtracting the sum up to $h-1$ from the sum up to $h$, we have
\begin{theorem}
In the limit $q\rightarrow\infty,$ and for $ h \leq n-4$,
\begin{equation}
|\sum_{\deg J=h}E_{\mu}(J,n,q)|= O(\frac{1}{q^{3/2}})
\end{equation}
\end{theorem}

\section*{Acknowledgements}
We gratefully acknowledge support under EPSRC Programme Grant EP/K034383/1
LMF: \textit{L}-Functions and Modular Forms.  JPK is also funded by a grant from the Leverhulme Trust, a Royal Society Wolfson Research Merit Award, and a Royal Society Leverhulme Senior Research Fellowship.   We are grateful to Professor Zeev Rudnick for helpful discussions and comments, and to a referee for a valuable suggestion.

\end{document}